\documentclass{article}
\usepackage{amssymb,latexsym,amsmath,amscd,amsthm,amsfonts}

\begin{document}

\begin{center}
{\large Automorphisms Inner in the Local Multiplier Algebra and Connes
Spectrum}

{\large \bigskip }

Raluca Dumitru, Costel Peligrad and Bogdan Visinescu

\textbf{\bigskip }
\end{center}

ABSTRACT: We prove that if $(A,G,\alpha)$ is a C*-dynamical system with $G$ abelian and if $g_{0}$ is in the inner centre of $G$ then $g_{0}\in \Gamma _{B}( \alpha ) ^{\perp }$. Conversely if $G$ is compact and $Sp( \alpha ) \diagup \Gamma ( \alpha )$ is finite, we prove that if $g_{0}\in \Gamma _{B}( \alpha)^{\perp }$ then $g_{0}$ is in the inner centre of $G$ (Proposition 3). We then prove that for $A$ prime and $G$ compact, not necessarily abelian, the
first part of Proposition 3 can be extended to the case when $\alpha_{g_{0}} $ is inner in $M_{loc}(A) ^{\alpha }$ (Proposition 7).

\bigskip


\bigskip

\bigskip AMS SUBJECT CLASSIFICATION: Primary 46L55, Secondary 46L40

\bigskip

\bigskip

\bigskip

\large In [8] and [4] , the C*-algebra of local multipliers was used for
the first time, under the name of algebra of essential (or limit)
multipliers. In [3], this algebra is extensively used in the problem of
"strong primeness" of fixed point algebras of C*-dynamical systems with
compact abelian groups. The recent book of P.Ara and M. Mathieu [1] contains
a systematic study of the algebra of local multipliers, $M_{loc}(A)$ and its applications.
\large

$M_{loc}(A)$ is by definition the inductive limit of the multiplier
algebras $M(I)$, where $I\subset A$ is an
essential ideal of $A$. If $I\subset J$ are two such
ideals and $z_{I}$, $z_{J}$ are the open central
projections in $A^{\prime \prime }$ corresponding to $I$,
respectively $J$ then the mapping $x\overset{\varphi
_{JI}}{\longrightarrow }${\large \ }$x\cdot z_{I}$ is an injective
morphism of $M(J)$ into $M(I)$. 
$M_{loc}(A)$ should be understood as the inductive limit of $
(M(I) ,\varphi _{JI})$.

In this paper we will obtain some results that relate the Connes spectrum with the innerness of automorphisms in the fixed point algebra $M_{loc}(A)^{\alpha }$.

The results are variants and/or extensions of corresponding results of Olesen, Pedersen and Stormer [6], [7], [9] and Ara and Mathieu [1] to the case of local multiplier algebras.

For definitions and basic results on Arveson spectrum, Connes spectrum and Borchers spectrum we send to [9].

In particular, if $(A,G,\alpha )$ is a C*-dynamical system with $G$ compact, abelian group and $\gamma \in \overset{\wedge }{G}$ we denote

\begin{equation*}
A_{\gamma } =\{x\in A\mid \alpha _{t}(x) =\langle t,\gamma \rangle x,\forall t\in G\}
\end{equation*}

For $\gamma=e$, the unit of $G$, we set $A_{e}=A^{\alpha}$.

\large Then the Arveson spectrum is defined by: \ 
\begin{equation*}
Sp( \alpha) =\{ \gamma \in \overset{\wedge }{G}\mid A_{\gamma }\neq ( 0)\}.
\end{equation*}

$\mathcal{H}^{\alpha }( A)$ denotes the set of all $\alpha$-invariant, hereditary C*-subalgebras of $A$
and $\mathcal{H}_{B}^{\alpha }(A)$ the set of all $\alpha$-invariant, hereditary subalgebras of $A$ ($C\in 
\mathcal{H}_{B}^{\alpha }(A)$ if $C\in \mathcal{H}^{\alpha }(A)$ and the ideal $\overline{ACA}$ is essential in A). Define the Connes spectrum and the Borchers spectrum as follows:

\begin{equation*}
\Gamma( \alpha)=\underset{C\in \mathcal{H}^{\alpha }(A)}{\bigcap}  Sp(\alpha \mid_{C})
\end{equation*}

\begin{equation*}
\Gamma_{B}(\alpha)=\underset{C\in \mathcal{H}_{B}^{\alpha }(A) }{\bigcap }Sp(\alpha \mid_{C})
\end{equation*}

Then $\Gamma(\alpha) \subseteq \Gamma _{B}(\alpha) \subseteq Sp(\alpha)$. $\Gamma(\alpha)$ is a subgroup of $\overset{\wedge }{G}$.

The corresponding spectra for actions of compact, non-abelian groups are defined in [5], [10].
\newtheorem{definition}{Definition}
\begin{definition}
If $(A,G,\alpha)$ is a C*-dynamical system we define $H\subset G$, the inner centre of G, as the set of all $t\in G$ for which there is an essential, $\alpha$-invariant ideal $I\subset A$ such that $\alpha _{t}=Ad$ $u$ for some $u\in M(I)^{\alpha}.$
\end{definition}
\newtheorem{remark}[definition]{Remark}
\begin{remark}
Since $M(I)^{\alpha }\subset M_{loc}(A)^{\alpha }$ for every essential ideal $I\subset A$, it is clear that every automorphism $\alpha _{t}$ with $t\in H$ is inner in $M_{loc}(A)^{\alpha }$. Here $M_{loc}(A)^{\alpha }$ denotes the fixed point algebra of $M_{loc}(A)$ under the extensions of all $\alpha_{t}$, $t\in G.$
\end{remark}

Our first result is concerned with automorphisms in the inner centre
of $G$.
\newtheorem{proposition}[definition]{Proposition}
\begin{proposition}
Let $(A,G,\alpha)$ be a C*-dynamical system with G locally compact abelian. If $t\in H$ (the inner centre of
G), then $\langle t,\sigma \rangle =1$ for every $\sigma \in \Gamma _{B}(\alpha)$. The converse holds if $G$ is compact and $Sp(\alpha) \diagup \Gamma(\alpha)$ is finite.
\end{proposition}

\begin{proof}
The proof of the direct implication is a modification of the corresponding result of Olesen ([9] Theorem 8.9.7) for multiplier algebras.

Let $\alpha _{t}=Ad$ $u$, $u\in M(I)^{\alpha }$ for some essential, $\alpha$-invariant ideal $I\subset A$. Then, since $\Gamma _{B}(\alpha) =\Gamma _{B}(\alpha \mid _{I})$ the proof of the first part of ([9] theorem 8.9.7) carries over word for word with $A$ replaced by $I$.

Conversely, let G be a compact group and $Sp(\alpha)\diagup \Gamma(\alpha)$ finite. We mention that we do not require $A$ to be simple as in ([9] Theorem 8.9.7).

Let $t_{0}\in G$ be such that $\langle t_{0},\sigma \rangle =1$, $\forall \sigma \in \Gamma _{B}(\alpha)$. Then, by ([9] Proposition 8.8.7) there exists $C\in \mathcal{H}_{B}^{\alpha}(A)$ such that $Sp(\alpha \mid _{C}) \subseteq \Gamma _{B}( \alpha ) +\Gamma ( \alpha)$.

Let $\sigma \in Sp(\alpha \mid _{C}).$ Then, there are $\sigma _{1}\in \Gamma _{B}( \alpha )$ and $\sigma _{2}\in \Gamma ( \alpha) \subseteq \Gamma_{B}( \alpha)$, such that $\sigma =\sigma _{1}+\sigma_{2}$. Then since $t_{0}$ annihilates $\Gamma _{B}(\alpha )$, it follows that $\langle t_{0},\sigma \rangle =1.$ Hence $\alpha _{t_{0}}\mid_{C}=id$. Let $I=\overline{ACA}$ be the two sided ideal generated by $C$. Since $C\in \mathcal{H}_{B}^{\alpha }(A)$, $I$ is an essential, $\alpha$-invariant ideal of $A$. Let $p\in C^{\prime \prime }\subset I^{\prime \prime }$ be the open projection associated to $C$. Then, the central support $c( p)$ in $I^{\prime \prime }$ equals 1. By ([9], Lemma 8.9.1) applied to $I^{\prime \prime }$ with $w=p$ we get a unique unitary $u\in I^{\prime \prime }$ such that $up=pu=p$ and $\alpha _{t_{0}}^{\prime \prime }\mid _{I^{\prime \prime }}=Ad$ $u$ and hence $\alpha_{t_{0}}\mid _{I}=Ad$ $u.$ 

Since $G$ is abelian, $\alpha _{t}^{\prime \prime}(u)$ satisfies all the properties of $u$. By uniqueness of $u$ with the above properties we get $\alpha
_{t}^{\prime \prime }(u) =u$ for every $t\in G$.

Define $I_{0}=\{ x\in I\mid ux \in I\}$. Then $I_{0}$ is a two sided, closed ideal of $I$. Indeed, $I_{0}$ is clearly a right ideal. It is also a left ideal since for $y \in I$ we have $uyx=uyu^{\ast }ux=\alpha
_{t_{0}}(y)ux\in I$. Since $p$ is the open projection associated with $C$, we have $px=x$ for every $x\in C$. Therefore, if $x\in C$, $ux=upx=px=x\in C\subset I$. Thus, $C\subseteq I_{0}$. Since the ideal generated by $C$ is $I$, it follows that $I_{0}=I$ and therefore $u$ is a left multiplier of $I$. Similarly, one can show that $u$ is a right multiplier of $I$ and so $u\in M(I)$. The proof is complete.
\end{proof}

{\large \bigskip }
\newtheorem{corollary}[definition]{Corollary}
\begin{corollary}
If $(A,G,\alpha)$ is a C*-dynamical system with $G$ compact abelian group and $Sp(\alpha)$ finite, then $\Gamma _{B}(\alpha) ^{\perp }$ $($the annihilator of $\Gamma _{B}( \alpha )$$)$ coincides with the inner centre. In particular, if $A$ is $G$-prime, then $\Gamma( \alpha) ^{\perp }( =\Gamma_{B}( \alpha) ^{\perp })$ coincides with the
inner centre.
\end{corollary}

This corollary generalizes ([1], Theorem 4.5.7 and the first part of
Theorem 4.5.8).

Our next result extends and improves on some results in [7], ([9]
Theorem 8.10.12).

\begin{proposition}
Let $(A,G,\alpha)$ be a C*-dynamical system with $G$ compact,abelian and $A$ prime. Suppose that $Sp(\alpha) \diagup \Gamma(\alpha)$ is finite. The following conditions are equivalent :

$(i)\ \ \ Sp( \alpha) =\Gamma(\alpha)$

$(ii)\ \ A^{\alpha }$ is prime

$(iii)\ M_{loc}( A) ^{\alpha }$ is prime

$(iv)$  The centre of $M_{loc}(A) ^{\alpha }$ is the scalars

$(v)$   No $\alpha _{t}\neq \iota $ is implemented by a unitary in $M_{loc}(A) ^{\alpha }$

$(vi)$  No $\alpha _{t}\neq \iota $ is in the inner centre of $G.$
\end{proposition}

\begin{proof}
$(i) \iff ( ii) $ by ([9] Theorem 8.10.4)

$( ii) \implies ( iii)$ : Since $A^{\alpha }$ is prime, every ideal of $A^{\alpha }$ is a prime C*-algebra. In particular, for every $\alpha$-invariant ideal $ I\subset A$, $I^{\alpha }$ is prime. Therefore, the
multiplier algebra $M(I^{\alpha })$ is prime. By ([3] Proposition 3.3), $M_{loc}( A)^{\alpha}=\lim_{I \text{invariant}   \; \text{ideal}}M(I^{\alpha })$. Since $M(I^{\alpha })$ are prime, their inductive limit, $M_{loc}(A) ^{\alpha }$ is prime (to show that inductive limits of prime algebras are prime, one can use [1], Lemma 1.2.32).

$(iii) \implies ( iv)$ : This is trivial.

$(iv) \implies ( v)$ : If $\alpha_{t_{0}}=Ad$ $u$ for $u\in M_{loc}(A)^{\alpha }$, then for $a\in M_{loc}(A)^{\alpha }$,  $\alpha _{t_{0}}(a) =a=uau^{\ast }$. Hence $u$ belongs to the centre of $M_{loc}(A)^{\alpha }$. Therefore $u=\lambda \cdot 1$ for some scalar and $\alpha _{t_{0}}=\iota $.

$(v) \implies (vi)$ : Evident since the inner centre is included in $M_{loc}(A)^{\alpha }.$

$(vi) \implies (i)$ : If $Sp(\alpha) \varsupsetneq \Gamma(\alpha)$ then since $\Gamma(\alpha)$ is a subgroup of $\overset{\wedge }{G}$, there exists $t_{0}\in \Gamma(\alpha)^{\bot }$ but $t_{0}\notin Sp(\alpha)^{\bot }$. By Proposition 3, $t_{0}$ belongs to the inner centre of $G$ i.e. there exists an $\alpha$-invariant ideal $I$ and a unitary $u\in M(I)^{\alpha }=M(I^{\alpha })$ such that $\alpha _{t_{0}}=Ad$ $u$. Therefore $\alpha _{t_{0}}=id$. Hence $t_{0}\in Sp(\alpha) ^{\bot}$, contradiction.
\end{proof}

{\large \bigskip }

Our next result is a stronger version of the first part of Proposition 3 for the case of a prime C*-algebra $A$ and a compact group $G$. Namely, we will prove that if $\alpha _{t_{0}}=Ad$ $u$, $u\in M_{loc}(A)^{\alpha }$ (rather than $u\in M(I)^{\alpha }$ for some essential, $\alpha$-invariant ideal) then $t_{0}\in \Gamma(\alpha)^{\bot } (=\Gamma _{B}(\alpha)^{\bot }$ since $A$ is now prime).

We will prove this for a compact, not necessarily abelian group $G$. The spectra $Sp(\alpha)$, $\Gamma(\alpha)$, $\overset{\sim }{\Gamma }(\alpha)$ for actions of compact, non abelian groups were defined in ([10], [5]). We will define the weak Connes spectrum $\Gamma _{\omega }(\alpha)$. Let $(A,G,\alpha)$ be a C*-dynamical system with $G$ compact. Let $\overset{\wedge }{G}$ be the dual of $G$ (i.e. the set of unitary equivalence classes of irreducible representations of $G$).

For $\pi \in \overset{\wedge }{G}$, $\pi =[ \pi _{ij}]$, $i,j=1,2,\ldots d_{\pi }$ ($d_{\pi }$ is the dimension of $\pi$) let $\chi _{\pi }(g) =\underset{i}{\sum }\pi _{ii}(g)$ be the character of $\pi$.

If $a\in A$, set $P(\pi)(a) =\underset{G}{\int }d_{\pi }\overline{\chi _{\pi }(g)} \alpha_{g}(a) dg$.

Denote
\begin{equation*}
A_{\pi }=\{ \underset{G}{\int }d_{\pi }\overline{\chi _{\pi }(g)} \alpha _{g}(a)dg\mid a\in A\}
\end{equation*}

Define 
\begin{equation*}
sp(\alpha ) =\{ \pi \in \overset{\wedge }{G}\mid A_{\pi}\neq (0) \}
\end{equation*}

The Connes spectrum $\Gamma(\alpha)$ and the strong Connes spectrum $\overset{\sim }{\Gamma }(\alpha)$ were defined in [5] and used in problems concerning the simplicity and primeness of the crossed products.

We now define the weak Connes spectrum:

\begin{equation*}
\Gamma_{\omega }(\alpha)=\underset{B\in \mathcal{H}^{\alpha }(A)}{\bigcap }sp(\alpha \mid _{B}).
\end{equation*}

{\large We have }$\overset{\sim }{\Gamma }\left( \alpha \right) \subseteq
\Gamma \left( \alpha \right) \subseteq \Gamma _{\omega }\left( \alpha
\right) .$

The first inclusion is strict even in the abelian case. The second inclusion is an equality for abelian groups but is strict for non-abelian groups (this is shown by the example 3.9 and Theorem 3.8 in [10] together with Theorem 2.5 in [5]).

Before starting our result, let us make some notations:

For $\pi =[\pi _{ij}] \in \overset{\wedge }{G}$, denote

\begin{equation*}
P_{ij}(\pi)(a)=\int  d_{\pi }\overline{\pi _{ji}(g)} \alpha _{g}(a)dg,\; a\in A.
\end{equation*}

Let $A_{ij}(\pi) =\{ P_{ij}(\pi)(a) \mid a\in A\}$.

We collect some facts about the bounded, linear mappings $P(\pi)$, $P_{ij}(\pi)$ from $A$ into $A$ in the following:

{\large \bigskip }
\newtheorem{lemma}[definition]{Lemma}
\begin{lemma}
$(i)$  $P(\pi) =\underset{i=\overline{1,d_{\pi }}}{\sum }P_{ii}(\pi)$

$(ii)$  $P_{ij}(\pi) P_{kl}(\pi) =\delta _{il}P_{kj}(\pi)$, where $\delta_{il}$ is the Kronecker symbol.

$(iii)$  $P(\pi)P_{ij}(\pi) =P_{ij}(\pi) =P_{ij}(\pi)P(\pi)$

$(iv)$  $P(\pi)^{2}=P(\pi)$

$(v)$  $\alpha _{g}(P_{ij}(\pi)( a) ) =\underset{k=\overline{1,d_{\pi }}}{\sum }P_{ik}(\pi)( a) \cdot \pi_{kj}(g)$

$(vi)$  If $P_{i_{0}j_{0}}(\pi)(a) \neq 0$, then $P_{i_{0}k}(\pi)(a) \neq 0$ for every $k=1,2,\ldots d_{\pi }$ and $ \{ P_{i_{0}k}(\pi )(a) \} _{k=1}^{d_{\pi}}$ are linearly independent.

$(vii)$  If $A_{\pi }\neq (0)$ then $A_{ij}(\pi ) \neq 0$ for every $i,j=1,2,\ldots d_{\pi }$.

$(viii)$  The subspaces $A_{\pi }\subset A$ are $\alpha$-invariant.
\end{lemma}

\begin{proof}
$(i)$ follows trivially from definitions.

$(ii)$, $( iii)$, $(iv)$ follow from $\pi( gg^{\prime }) =\pi(g) \pi(g^{\prime }) $ ($\pi$ is a representation of $G$) and from Weil's orthogonality relations ([2], 2.2.39).

$(v)$ is easily seen since $\pi$ is a representation, $\pi(gg^{\prime }) =\pi(g) \pi( g^{\prime })$ ([2], 2.2.39).

$(vi)$ If $P_{i_{0}k}( \pi)(a) =0$, then by $(ii)$, $P_{kj_{0}}(\pi)(a) P_{i_{0}k}(\pi)(a) =P_{i_{0}j_{0}}(\pi)(a)=0$, contradiction. So $P_{i_{0}k}(\pi)(a) \neq 0$, $k=1,2,\ldots d_{\pi }$ if $P_{i_{0}j_{0}}(\pi)(a) \neq 0$. Let now $\underset{k}{\sum }\lambda _{k}P_{i_{0}k}(\pi)(a) =0$. Applying $P_{k_{0}j_{0}}(\pi)$ to the previous equality, it follows that $\lambda _{k_{0}}P_{i_{0}j_{0}}(\pi)(a) =0$ and hence $\lambda _{k_{0}}=0$ for every $k_{0}=1,2,\ldots d_{\pi }$.

$(vii)$ If $A_{\pi }\neq 0$, then there is $a\in A$ such that $P(\pi)(a) =\sum P_{ii}(\pi)(a) \neq 0$. Hence $P_{i_{0}i_{0}}(\pi)(a) \neq 0$ for some $i_{0}$.

{\large By }$\left( vi\right) ${\large \ }$P_{i_{0}k}\left( \pi \right)
\left( a\right) \neq 0${\large \ for every }$k=1,2,\ldots d_{\pi }.$

Let $b=P_{i_{0}i}(\pi) ( a)$. Then $b\neq 0$ and $P_{ij}( \pi)( b)=P_{i_{0}j}(\pi)( a) \neq 0$. So $P_{ij}(\pi)(b) \in A_{ij}(\pi) \neq ( 0)$.

$(viii)$ Follows from $(ii)$, $(iii)$, $(iv)$,  and $(v)$.
\end{proof}

We can now prove:

\begin{proposition}
Let $(A,G,\alpha)$ be a C*-dynamical system with $G$ compact and $A$ prime. If $\alpha _{g_{0}}=Ad$ $u_{0}$ for some $u_{0}\in M_{loc}(A) ^{\alpha }$, then $g_{0}\in \Gamma _{\omega}(\alpha)^{\bot }$ $(\subseteq \Gamma(\alpha)^{\bot} \subseteq \overset{\sim }{\Gamma }(\alpha)^{\bot }$).
\end{proposition}

\begin{proof}
Let $\epsilon >0$. By ([3] Proposition 3.3 and [11] Proposition L.2.2) there exists an $\alpha$-invariant ideal $I\subset A$ and an unitary $u\in M(I)^{\alpha }$ such that $\Vert u_{0}-u\Vert <\frac{\epsilon }{4}$.
Let $\lambda _{0}\in sp( u)$ and $f$ be a positive continuous function on the unit circle $\mathbb{T}$ with $supp$ $f\subseteq \{ \lambda \in \mathbb{T}\mid \vert \lambda -\lambda _{0}\vert <\frac{\epsilon }{4}\}$. Let also $h$ be a positive, continuous function on the unit circle such that $0\leq h\leq 1$, $h\mid _{\text{supp}f}\equiv 1$ and $h(\lambda) =0$ for every $\lambda \in \mathbb{T}$ with $\vert \lambda -\lambda _{0}\vert \geq \frac{\epsilon }{4}$. Then $f(u) ,h(u) \in M( I) ^{\alpha }$ and $f( u) h(u) =f(u)$. Let $C=\overline{f(u)If(u)}$ (norm closure). Then $C\in \mathcal{H}^{\alpha }(I) \subseteq \mathcal{H}^{\alpha}(A)$.

Let $x\in C$. Then:
\begin{eqnarray*}
\Vert \alpha _{g_{0}}(x)-x\Vert&=&\Vert u_{0}xu_{0}^{\ast }-x\Vert =\Vert u_{0}x-xu_{0}\Vert \\
&\leq &\Vert ( u_{0}-u) x-x( u_{0}-u) \Vert+\Vert ux-xu\Vert \\
&<&\frac{\epsilon }{2}\Vert x\Vert +\Vert ux-xu\Vert
\end{eqnarray*}

Since $x\in C=\overline{f( u)If( u) }=\overline{h( u) f(u)If(u)h(u) }$ we have: 
\begin{equation*}
\Vert ux-xu\Vert =\Vert ( u-\lambda _{0})x-x(u-\lambda _{0}) \Vert <\frac{\epsilon }{2}\Vert x\Vert
\end{equation*}

Therefore:
\begin{equation}
\Vert \alpha _{g_{0}}( x) -x\Vert <\epsilon \Vert x\Vert
\end{equation}

Let now $\pi \in \Gamma _{\omega }( \alpha )$. Since $C\in \mathcal{H}^{\alpha }(A)$, it follows that $\pi \in Sp(\alpha)$. Then $C_{\pi }\neq (0)$. By Lemma 6 $(vii)$ there is $a\in C$ such that $P_{ii}(\pi)(a) =a_{ii}\neq 0$ for every $i=1,2,\ldots d_{\pi }$.

From relation (1) with $x=a_{ii}$, we have:
\begin{equation}
\Vert \alpha _{g_{0}}( a_{ii}) -a_{ii}\Vert <\epsilon \Vert a_{ii}\Vert
\end{equation}

Lemma 6 $(v)$ and (2) above give:
\begin{equation}
\Vert \underset{k=\overline{1,d_{\pi }}}{\sum }a_{ik}\pi _{ki}(g_{0}) -a_{ii}\Vert <\epsilon \Vert a_{ii}\Vert \label{a}
\end{equation}

where $a_{ik}=P_{ik}(\pi) ( a_{ii})$.

By Lemma 6 $(vi)$ and Hahn-Banach there is a functional $\varphi$ of norm one on $C$ such that $\varphi(a_{ii})=\Vert a_{ii}\Vert$ and $\varphi( a_{ik}) =0$ for all $k\neq i$.

Applying $\varphi $ to (3), we get:
\begin{equation}
\vert \pi _{ii}( g_{0}) -1\vert \cdot \Vert a_{ii}\Vert <\epsilon \Vert a_{ii}\Vert
\end{equation}

Hence $\vert \pi _{ii}( g_{0}) -1\vert <\epsilon$ for every $\epsilon >0$. Therefore $\pi_{ii}(g_{0})=1$ for all $i=1,2,\ldots d_{\pi }$.
Since $\pi ( g_{0}) =[ \pi _{ij}( g_{0})]$ is a unitary matrix with all diagonal entries equal 1, it follows immediately that $\pi( g_{0}) =1$ and we are
done.
\end{proof}

\newpage
\textbf{Acknowledgements}. R. Dumitru was supported by "Isabel and
Mary Neff" Fellowship and by URC Summer Grant. C. Peligrad was supported by
a Taft Travel Grant. B. Visinescu was supported by a Taft Enhancement
Fellowship.

{\large \bigskip }

{\large \bigskip }

$\bigskip $

$%
\begin{array}{ccc}
\text{Costel Peligrad} &  & \text{ Raluca Dumitru and Bogdan
Visinescu} \\ 
Department\text{ }of\text{ }Mathematics &  & Department\text{ }of\text{ }%
Mathematics \\ 
University\text{ }of\text{ }Cincinnati &  & University\text{ }of\text{ }%
Cincinnati \\ 
Cincinnati,{\large \ }Ohio,{\large \ }45221-0025 &  & Cincinnati,{\large \ }%
Ohio,{\large \ }45221-0025 \\ 
E-mail:{\large \ }peligrc@math.uc.edu &  & E-mail:{\large \ }%
dumitrur@math.uc.edu \\ 
&  & E-mail:{\large \ }visinebc@math.uc.edu \\ 
&  &  \\ 
&  & On\text{ }leave\ from: \\ 
&  & Institute\text{ }of\text{ }Mathematics \\ 
&  & of\text{ }the\text{ }Romanian\text{ }Academy \\ 
&  & P.O.\ Box\ 1-764 \\ 
&  & 70700\ Bucharest \\ 
&  & Romania%
\end{array}%
$

{\large \bigskip }

{\large \bigskip }

{\large \bigskip }

\end{document}